\documentclass{elsart}
\usepackage{natbib}
\usepackage{amsfonts,amsmath,amssymb}
\usepackage{dsfont}
\usepackage{graphicx}

\newenvironment{disarray}
{\everymath{\displaystyle\everymath{}}\array}
{\endarray}

\def\findemo{\hfill$\Box$}
\def\E{\mathbb{E}}
\def\pen{{\mathrm{pen}}}
\def\1{\mathds{1}}

\begin{document}

\begin{frontmatter}

\title{Nonparametric estimation of the stationary density and the transition density 
of a Markov chain}

\author{Claire Lacour}

\address{ MAP5, Universit\'e Paris Descartes, CNRS UMR 8145, 45 rue des 
          Saints-P\`eres 75270 Paris Cedex 06, France \\
lacour@math-info.univ-paris5.fr}

\begin{abstract}
In this paper, we study first the problem of nonparametric estimation of the stationary density 
$f$ of a discrete-time Markov chain $(X_i)$. We consider a collection of projection estimators on 
finite dimensional linear spaces. We select an estimator among the collection by minimizing 
a penalized contrast. The same technique enables to estimate the density $g$ of $(X_i, X_{i+1})$ 
and so to provide an adaptive estimator of the transition density $\pi=g/f$.
We give bounds in $L^2$ norm for these estimators and we show that they are
adaptive in the minimax sense over a large class of Besov spaces. Some examples and simulations 
are also provided.
\end{abstract}


\end{frontmatter}

\section{Introduction}
\label{intro}

Nonparametric estimation is now a very rich branch of statistical theory. The case of i.i.d. 
observations is the most detailed but many authors are also interested in the case of Markov 
processes. Early results are stated by \citet{Rou69b}, who studies nonparametric estimators of the 
stationary density and the transition density of a Markov chain. He considers kernel estimators 
and assumes that the chain satisfies the strong Doeblin's condition $(D_0)$ (see \citet{doob} p.221). 
He shows consistency and asymptotic normality of his estimator. Several authors tried 
to consider weaker assumptions than the Doeblin's condition. \citet{rosen70} introduces an other 
condition, denoted by $(G_2)$, and he gives results on the bias and the variance of the kernel 
estimator of the invariant density in this weaker framework. \citet{yako89} improves also the result 
of asymptotic normality by considering a Harris-condition. The study of kernel estimators is 
completed by \citet{masry&gyorfi} who find sharp rates for this kind of estimators of the 
stationary density and by \citet{basu&sahoo} who prove a Berry-Esseen inequality under the condition 
$(G_2)$ of Rosenblatt. Other authors are interested in the estimation of the invariant distribution 
and the transition density in the non-stationary case: \citet{dou&ghin} bound the integrated risks 
for any initial distribution. In \citet{hernandez}, recursive estimators for a 
non-stationary Markov chain are described. \citet{liebscher} gives results for the invariant density 
in this non-stationary framework using a  condition denoted by $(D_1)$ derived from the Doeblin's 
condition but weaker than $(D_0)$. All the above papers deal with kernel estimators. Among those 
who are not interested in such estimators, let us mention \citet{bosq} who studies an estimator of 
the stationary density by projection on a Fourier basis, \citet{rao} who outlines a new estimator 
for the stationary density by using delta-sequences and \citet{gillertwartenberg} who present 
estimators based on Hermite bases or trigonometric bases.

The recent work of \citet{theseclem} allows to measure the performance of all these estimators 
since he proves lower bounds for the minimax rates and gives thus the optimal convergence rates
for the estimation of the stationary density and the transition density.
Cl\'emen\c{c}on also provides an other kind of estimator for the stationary density
and for the transition density, that he 
obtains by projection on wavelet bases. He presents an adaptive procedure which is 
"quasi-optimal" in the sense that the procedure reaches almost the optimal rate but with a 
logarithmic loss. He needs other conditions than those we cited above and in particular a 
minoration condition derived from Nummelin's (1984) works. In this paper, we will use the same condition. 

The aim of this paper is to estimate the stationary density of a discrete-time Markov chain and its 
transition density. 
We consider an irreducible positive recurrent Markov chain $(X_n)$ with a stationary density denoted 
by $f$. We suppose that the initial density is $f$ (hence the process is stationary) and we 
construct an estimator $\tilde{f}$ from the data $X_1,\dots, X_n$. Then, we study the mean 
integrated squared error $\E\|\tilde{f}-{f}\|_2^2$ and its convergence rate. 
The same technique enables to estimate the density $g$ of $(X_i, X_{i+1})$ 
and so to provide an estimator of the transition density $\pi=g/f$, called the quotient estimator.

An adaptative procedure is proposed for the two estimations and it is proved that both resulting 
estimators reach the optimal minimax rates without additive logarithmic factor.

We will use here some technical methods known as the Nummelin splitting technique (see \citet{numm},
\citet{m&t} or \citet{loscher}). This method allows to reduce the general state space Markov chain
theory to the countable space theory. Actually, the splitting of the original chain creates an 
artificial accessible atom and we will use the hitting times to this atom to decompose the chain, 
as we would have done for a countable space chain.

To build our estimator of $f$, we use model selection via penalization as described in \citet{BBM}.
First, estimators by projection denoted by $\hat{f}_m$ are considered. The index $m$ denotes the 
model, i.e. the subspace to which the estimator belongs. Then the model selection technique allows 
to select automatically an estimator $\hat{f}_{\hat{m}} $ from the collection 
of estimators $(\hat{f}_m)$. The estimator of $g$ is built in the same way.
The collections  of models that we consider here include wavelets but 
also trigonometric polynomials and piecewise polynomials.

This paper is organized as follows. In Section 2, we present our assumptions on the Markov 
chain and on the collections of models. We give also examples of chains and models. 
Section 3 is devoted to estimation of the stationary density and in Section 4 the estimation of the 
transition density is explained. Some simulations are presented in Section 5.
The proofs are gathered in the last section, which contains also a presentation of the 
Nummelin splitting technique.

\section{The framework}

\subsection{Assumptions on the Markov chain}
We consider an irreducible Markov chain $(X_n)$ taking its values in the real 
line $\mathbb{R}$. We suppose that $(X_n)$ is positive recurrent, i.e. it admits a 
stationary probability measure $\mu$ (for more details, we refer to \citet{m&t}). We assume that the 
distribution $\mu$ has a density $f$ with respect to the Lebesgue measure and it is this 
quantity that we want to estimate. Since the number of observations is finite, $f$ is estimated 
on a compact set only. Without loss of generality, this compact set is assumed to be equal to 
$[0,1]$ and, from now, $f$ denotes the transition density multiplied by the indicator 
function of $[0,1]$ $f\1_{[0,1]}$. 
More precisely, the Markov process is supposed to satisfy the following assumptions:
\begin{enumerate}
\item[A1.] $(X_n)$ is irreducible and positive recurrent.
\item[A2.] The distribution of $X_0$ is equal to $\mu$ , thus the chain is (strictly) stationary.
\item[A3.] The stationary density $f$ belongs to $L^\infty([0,1])$ i.e.
$\sup_{x\in[0,1]}|f(x)|<\infty$
\item[A4.] The chain is strongly aperiodic, i.e. it satisfies the following minorization condition: 
there is some function $h : [0,1]\mapsto [0,1]$ with $\int hd\mu>0$ and a positive distribution 
$\nu$ such that, for all event $A$ and for all $x$,
$$P(x,A)\geq h(x)\nu(A)$$
where $P$ is the transition kernel of $(X_n)$.
\item[A5.] The chain is geometrically ergodic, i.e. there exists a function $V>0$ finite 
and a constant $\rho\in (0,1)$ such that, for all $n\geq 1$
$$\|P^n(x,.)-\mu\|_{TV}\leq V(x)\rho^n$$ where $\|.\|_{TV}$ is the total variation norm.
\end{enumerate}

We can remark that condition A3 implies that $f$ belongs to $L^2([0,1])$ where
$\displaystyle L^2([0,1])=\{t:\mathbb{R}\mapsto\mathbb{R}, \text{Supp}(t)\subset[0,1] \text{ and  } 
\|t\|^2=\int_0^1 t^2(x)dx<\infty\}.$

Notice that, if the chain is aperiodic, condition A4 holds, at least for some $m$-skeleton 
(i.e. a chain with transition probability $P^m$) (see Theorem 5.2.2 in \citet{m&t}).
This minorization condition is used in the Nummelin splitting technique and is also required in 
\citet{theseclem}.

The last assumption, which is called geometric regularity by \citet{clem}, means that the 
convergence of the chain to the invariant distribution is geometrically fast.
In \citet{m&t}, we find a slightly 
different condition (replacing the total variation norm by the $V$-norm). This condition, which
is sufficient for A5, 
is widely used in Monte Carlo Markov Chain literature because it guarantees central 
limit theorems and enables to simulate laws via a Markov chain (see for example \citet{jarner},
\citet{robertsrosen98} or \citet{meyn94}).

The following subsection gives some examples of Markov chains satisfying hypotheses A1--A5.

\subsection{Examples of chains}

\subsubsection{Diffusion processes}
We consider the process $(X_{i\Delta })_{1\leq i \leq n}$ where 
$\Delta>0$ is the observation step and $(X_t)_{t\geq 0}$ is defined by
$$dX_t=b(X_t)dt +\sigma(X_t)dW_t$$
where $W$ is the standard Brownian motion, $b$ a
is a locally bounded Borelian function and $\sigma$ is a uniformly continuous function  such that:
\begin{enumerate}
\item there exists $\lambda_-$, $\lambda_+$ such that $\forall x\neq 0$, $0<\lambda_-<\sigma^2(x)<\lambda_+$,
\item there exists $M_0,\alpha\geq 0$ and $r>0$ such that 
$\forall |x|\geq M_0, xb(x)\leq -r|x|^{\alpha+1}.$
\end{enumerate}

Then, if $X_0$ follows the stationary distribution, Proposition 1 in \cite{pardouxveretennikov01}
shows that the discretized process 
$(X_{i\Delta })_{1\leq i \leq n}$ satisfies Assumptions A1--A5.

\subsubsection{Nonlinear AR(1) processes}
Let us consider the following process
$$X_n=\varphi(X_{n-1})+\varepsilon_{X_{n-1},n}$$
where $\varepsilon_{x,n}$ has a positive density  $l_x$ with respect to the Lebesgue measure, 
which does not depend on $n$. We suppose that $\varphi$ is bounded on any compact set 
and that there exist $M>0$ and $\rho<1$ such that, for all $|x|>M$, $|\varphi(x)|<\rho |x|$.
\citet{Mokkadem} proves that if there exists $s>0$ such that $\sup_x\E|\varepsilon_{x,n}|^s<\infty$,
then the chain is geometrically ergodic.
If we assume furthermore that $l_x$ has a lower bound then the chain satisfies all the previous
assumptions.

\subsubsection{ARX (1,1) models}
The nonlinear process ARX(1,1) is defined by
$$X_n=F(X_{n-1},Z_n) + \xi_n$$
where $F$ is bounded and $(\xi_n)$, $(Z_n)$ are independent sequences of i.i.d.
random variables with $\E|\xi_n|<\infty$. We suppose that the distribution of  $Z_n$ has a positive  
density $l$ with respect to the Lebesgue mesure.
Assume that there exist $\rho<1$, a locally bounded and mesurable function 
$h:\mathbb{R}\mapsto\mathbb{R}^+$ such that $\E h(Z_n)<\infty$ and positive constants $M, c$ 
such that 
$$\forall |(u,v)|>M \quad |F(u,v)|<\rho |u|+h(v)-c \text{  and }
\sup_{|x|\leq M}|F(x)|<\infty.$$
Then the process $(X_n)$ satisfies Assumptions A1--A5 (see \citet{doukhan} p.102).

\subsubsection{ARCH process}
The considered model is
$$X_{n+1}=F(X_n)+ G(X_n)\varepsilon_{n+1}$$
where $F$ and $G$ are continuous functions and for all $x$, $G(x)\neq 0$.
We suppose that the distribution of  $\varepsilon_n$ has a positive and continuous density 
with respect to the Lebesgue measure 
and that there exists $s\geq 1$ such that $\E|\varepsilon_n|^s<\infty$.
The chain $(X_i)$ satisfies Assumptions A1--A5 if (see \citet{doukhan} p.106): 
$$\lim\sup_{|x|\to\infty} \frac{|F(x)|+|G(x)|(\E|\varepsilon_n|^s)^{1/s}}{|x|}<1.$$

\subsection{Assumptions on the models}\label{hypmo}
In order to estimate $f$, we need to introduce some collections of models.
The assumptions on the models are the following:

\begin{enumerate}
\item[M1.] Each $S_m$ is a linear subspace of $(L^\infty\cap L^2)([0,1])$ with dimension 
$D_m\leq\sqrt{n}$
\item[M2.] Let  $$\phi_m=\frac{1}{\sqrt{D_m}}\underset{t\in S_m\backslash\{0\}}
{\sup}\frac{\|t\|_\infty}{\|t\|}$$
There exists a real $r_0$ such that for all $m,\quad\phi_m\leq r_0.$
\end{enumerate}

This assumption ($L^2$-$L^\infty$ connexion) is introduced by \citet{BBM} and can be written:
\begin{equation}
\forall t\in S_m \qquad \|t\|_\infty\leq r_0\sqrt{D_m}\|t\|.\label{M2}
\end{equation}

\noindent We get then a set of models $(S_m)_{m\in\mathcal{M}_n}$ where 
$\mathcal{M}_n=\{m,\quad D_m\leq\sqrt{n}\}$. We need now a last assumption regarding the whole
collection, which ensures that, for $m$ and $m'$ in $\mathcal{M}_n$,
$S_m+S_m'$ belongs to the collection of models.

\begin{enumerate}
\item[M3.]The models are nested, that is for all $m$, $D_m\leq D_{m'}\Rightarrow S_m\subset S_{m'}$.
\end{enumerate}

\subsection{Examples of models}\label{ex}
We show here that the assumptions M1-M3 are not too restrictive. Indeed, they are verified for 
the models spanned by the following bases (see \citet{BBM}):

\begin{itemize}

\item Histogram basis:
 $ S_m=<\varphi_1,\dots,\varphi_{2^m}>$ with $\varphi_j
=2^{m/2}\1_{[\frac{j-1}{2^m},\frac{j}{2^m}[}  $ for $j=1,\dots,2^m$. Here $D_m=2^m$, $r_0=1$ and
$\mathcal{M}_n=\{1,\dots,\lfloor{\ln n}/{2\ln2}\rfloor\}$ where $\lfloor x\rfloor$ denotes the 
floor of $x$, i.e. the largest integer less than or equal to $x$.

\item Trigonometric basis: 
$S_m=<\varphi_0,\dots,\varphi_{m-1}>$ with $\varphi_0(x)=\1_{[0,1]}(x)$, 
$\varphi_{2j}=\sqrt{2}$ $\cos(2\pi jx)\1_{[0,1]}(x)$, $\varphi_{2j-1}=\sqrt{2}\sin(2\pi jx)
\1_{[0,1]}(x)$ for $j\geq 1$. For this model $D_m=m$ and $r_0=\sqrt{2}$ hold.

\item Regular piecewise polynomial basis: 
$S_m$ is spanned by polynomials of degree $0,\dots,r$ (where $r$ is fixed) on each interval 
$[(j-1)/2^D,j/2^D[, j=1,\dots,2^D$. In this case, $m=(D,r)$, $D_m=(r+1)2^D$ and 
$\mathcal{M}_n=\{(D,r),~ D=1,\dots,\lfloor\log_2(\sqrt{n}/(r+1)) \rfloor\}$.We can put 
$r_0=\sqrt{r+1}$.

\item Regular wavelet basis:
$S_m=<\psi_{jk}, j=-1,\dots,m, k\in\Lambda(j)>$ where $\psi_{-1,k}$ points out the translates of 
the father wavelet and $\psi_{jk}(x)=2^{j/2}\psi(2^jx-k)$ where $\psi$ is the mother wavelet.
We assume that the support of the wavelets is included in $[0,1]$ and that $\psi_{-1}=\varphi$ 
belongs to the Sobolev space $W_2^r$.
In this framework $\Lambda(j)=\{0,\dots,K2^j-1\}$ (for $j\geq 0$) where $K$ is a constant which 
depends on the supports of $\varphi$ and $\psi$: for example for the Haar basis $K=1$.
We have then  $D_m=\sum_{j=-1}^m |\Lambda(j)|=|\Lambda(-1)|+K(2^{m+1}-1)$. Moreover
\begin{eqnarray*}
\phi_m &\leq &\frac{\sum_k|\psi_{-1,k}|+\sum_{j=0}^m2^{j/2} \sum_k|\psi_{j,k}|}{\sqrt{D_m}}\\
   &\leq& \frac{\|\varphi\|_\infty\vee\|\psi\|_\infty(1+\sum_{j=0}^m2^{j/2})}
          {\sqrt{(K\wedge|\Lambda(-1)|)2^{m+1}}} 
   \leq \frac{\|\varphi\|_\infty\vee\|\psi\|_\infty}{K\wedge|\Lambda(-1)|}=:r_0
\end{eqnarray*}

\end{itemize}

\section {Estimation of the stationary density }

\subsection{Decomposition of the risk for the projection estimator }

Let \begin{equation}
\gamma_{n}(t)=\frac{1}{n}\sum_{i=1}^{n}[\|t\|^2-2t(X_i)].\label{gamman}
\end{equation}
Notice that $\E(\gamma_{n}(t))=\|t-f\|^2-\|f\|^2$ and therefore $\gamma_{n}(t)$ is 
the empirical version of the $L^2$ distance between $t$ and $f$.
Thus, $\hat{f}_m$ is defined by 
\begin{equation}
\hat{f}_m=\underset{t\in S_m}{\arg\min} \gamma_{n}(t)\label{defest}
\end{equation}
where $S_m$ is a subspace of $L^2$ which satisfies M2. 
Although this estimator depends on $n$, no index $n$ is mentioned in order 
to simplify the notations . It is also the case for all the estimators in this paper.

A more explicit formula for $\hat{f}_m$ is easy to derive:
\begin{equation}
\hat{f}_m=\sum_{\lambda\in\Lambda}\hat{\beta}_\lambda\varphi_\lambda,
\qquad \hat{\beta}_\lambda=\frac{1}{n}\sum_{i=1}^{n}\varphi_\lambda(X_i)
\label{beta}\end{equation}
where $(\varphi_\lambda )_{\lambda\in\Lambda}$ is an orthonormal basis of $S_m$.
Note that $$\E(\hat{f}_m)=\sum_{\lambda\in\Lambda}<f,\varphi_\lambda>\varphi_\lambda,$$ 
which is the projection of $f$ on $S_m$.

In order to evaluate the quality of this estimator, we now compute the mean 
integrated squared error $\E\|f-\hat{f}_m\|^2$ (often denoted by MISE).

\begin{prop} \label{est1}
Let $X_n$ be a Markov chain which satisfies Assumptions A1--A5 and $S_m$ be a subspace
of $L^2$ with dimension $D_m\leq n$. If $S_m$ satisfies condition M2, then the estimator 
$\hat{f}_m$ defined by \eqref{defest} satisfies
$$\E\|f-\hat{f}_m\|^2\leq d^2(f,S_m)+C\frac{D_m}{n}$$
where $C$ is a constant which does not depend on $n$.
\end{prop}

To compute the bias term $d(f,S_m)$, we assume that $f$ belongs to the Besov space 
$B_{2,\infty}^\alpha([0,1])$. We refer to \cite{devorelorentz} p.54 for the definition of $B_{2,\infty}^\alpha([0,1])$.
Notice that when $\alpha$ is an integer, the Besov space $B_{2,\infty}^\alpha([0,1])$ contains the 
Sobolev space $W_2^\alpha$ (see \citet{devorelorentz} p.51--55).

Hence, we have the following corollary.

\begin{cor}\label{coro1}
Let $X_n$ be a Markov chain which satisfies Assumptions A1--A5. 
Assume that the stationary density $f$ belongs to 
$B_{2,\infty}^\alpha([0,1])$ and that $S_m$ is one of the spaces mentioned in
Section \ref{ex} (with the regularity of polynomials and wavelets larger than $\alpha-1$). 
If we choose $D_m=\lfloor n^\frac{1}{2\alpha+1}\rfloor$, then the estimator defined by 
\eqref{defest} satisfies
$$\E\|f-\hat{f}_m\|^2=O(n^{-\frac{2\alpha}{2\alpha+1}})$$
\end{cor}

We can notice that we obtain the same rate than in the i.i.d. case (see \citet{DJKP}). Actually, 
\citet{theseclem} proves that $n^{-\frac{2\alpha}{2\alpha+1}}$ is the optimal rate in the minimax 
sense in the Markovian framework. 
With very different theoretical tools, \citet{tribouleyviennet} show that this rate is also reached 
in the case of the univariate density estimation of $\beta$-mixing random variables by using a 
wavelet estimator.

However, the choice $D_m=\lfloor n^\frac{1}{2\alpha+1}\rfloor$ is possible only 
if we know the regularity $\alpha$ of the unknown $f$. But generally, it is not the case. It is 
the reason why we construct an adaptive estimator, i.e. an estimator which achieves the optimal 
rate without requiring the knowledge of $\alpha$. 

\subsection {Adaptive estimation }

Let $(S_m)_{m\in\mathcal{M}_n}$ be a collection of models as described in Section \ref{hypmo}. 
For each $S_m$, $\hat{f}_m$ is defined as above by \eqref{defest}.
Next, we choose $\hat{m}$ among the family $\mathcal{M}_n$ such that
$$\hat{m}=\underset{m\in \mathcal{M}_n}{\arg\min} [\gamma_{n}(\hat{f}_m)+\pen(m)] $$
where $\pen$ is a penalty function to be specified later.
We denote $\tilde{f}=\hat{f}_{\hat{m}}$ and we bound the $L^2$-risk $\E\|f-\tilde{f}\|$ as follows.

\begin{thm}\label{th1}
Let $X_n$ be a Markov chain which satisfies Assumptions A1--A5 and $(S_m)_{m\in\mathcal{M}_n}$ 
be a collection of models satisfying Assumptions M1--M3. Then the estimator defined by 
\begin{equation} 
\tilde{f}=\hat{f}_{\hat{m}} \quad\text{ where } \quad
\hat{m}=\underset{m\in \mathcal{M}_n}{\arg\min} [\gamma_{n}(\hat{f}_m)\label{tildef}+\pen(m)],
\end{equation} with 
\begin{equation}\pen(m)=K\frac{D_m}{n} \text{ for some  } K> K_0\label{penalite}\end{equation} 
(where $K_0$ is a constant depending on the chain)
satisfies
$$\mathbb{E}\|\tilde{f}-f\|^2 \leq 3\underset{m\in\mathcal{M}_n}{\inf}\{d^2(f,S_m) 
    +\pen(m)\}+\frac{C_1}{n}$$
where $C_1$ does not depend on $n$.
\end{thm}
 
\begin{rem}\label{rempenalite}
The constant $K_0$ in the penalty depends only on the distribution of the chain and can be chosen
equal to $\max(r_0^2,1)(C_1+C_2\|f\|_\infty)$ where $C_1$ and $C_2$ are 
theoretical constants provided by the Nummelin splitting technique.
The number $r_0$ is known and depends on the chosen base (see subsection \ref{hypmo}).
The mention of $\|f\|_\infty$ in the penalty term seems to be a problem, seeing that $f$ is unknown.
Actually, we could replace $\|f\|_\infty$ by $\|\hat{f}\|_\infty$ with $\hat{f}$ an estimator of $f$.
This method of random penalty is successfully applied in \cite{birge&massart97} or \cite{comte2001} 
for example. But we choose not to use this method here, since 
the constants $C_1$ and $C_2$ in $K_0$ are not computable either. 
Notice that \cite{clem} handle with the same kind of unknown quantities in the threshold of his nonlinear
wavelet estimator. Actually it is the price to pay for dealing with dependent variables 
(see also the mixing constant in the threshold in \cite{tribouleyviennet}).
But this annoyance can be circumvented for practical purposes.
Indeed, for the simulations the computation of the penalty is hand-adjusted. Some techniques 
of calibration can be found in \cite{lebarbier05} in the context of multiple change point detection.
In a Gaussian framework the practical choice of the penalty for implementation 
is also discussed in Section 4 of \cite{birgemassart05}.
\end{rem}

\begin{cor}\label{coro2}
Let $X_n$ be a Markov chain which satisfies Assumptions A1--A5 and $(S_m)_{m\in\mathcal{M}_n}$ 
be a collection of models mentioned in
Section \ref{ex} (with the regularity of polynomials and wavelets larger than $\alpha-1$).
If $f$ belongs to 
$B_{2,\infty}^\alpha([0,1])$, with $\alpha>1/2$, then the estimator defined by 
\eqref{tildef} and \eqref{penalite}
satisfies $$\mathbb{E}\|\tilde{f}-f\|^2 =O(n^{-\frac{2\alpha}{2\alpha+1}})$$
\end{cor}

\begin{rem}
When $\alpha>\frac{1}{2}$, $B_{2,\infty}^\alpha([0,1])\subset C[0,1]$ (where $C[0,1]$ is the set of 
the continuous functions with support in $[0,1]$) and then the assumption A3 $\|f\|_\infty <\infty$  
is superfluous.
\end{rem}
  
We have already noticed that it is the optimal rate in the minimax sense (see the lower bound in 
\citet{theseclem}).
Note that here the procedure reaches this rate whatever the regularity of $f$, without needing to 
know $\alpha$. This result is thus a improvement of the one of \citet{theseclem}, whose 
adaptive procedure achieves only the rate $(\log(n)/n)^{\frac{2\alpha}{2\alpha+1}}$. Moreover, our 
procedure allows to use more bases (not only wavelets) and is easy to implement.

\section{Estimation of the transition density}

We now suppose that the transition kernel $P$ has a density $\pi$.
In order to estimate $\pi$, we remark that $\pi$ can be written $g/f$ where
$g$ is the density of $(X_i, X_{i+1})$.
Thus we begin with the estimation of $g$. 
As previously, $g$ and $\pi$ are estimated on a compact set which is assumed to be equal to 
$[0,1]^2$, without loss of generality.

\subsection{Estimation of the joint density $g$}

We need now a new assumption.
\begin{enumerate}
\item[A3'.]  $\pi$ belongs to $L^\infty([0,1]^2)$.
\end{enumerate}
Notice that A3' implies A3. We consider now the following subspaces.
$$S_m^{(2)}=\{t\in L^2([0,1]^2),\quad t(x,y)=\sum_{\lambda,\mu\in\Lambda_m}
\alpha_{\lambda,\mu}\varphi_\lambda(x)\varphi_\mu(y)\}$$ where
$(\varphi_\lambda)_{\lambda\in\Lambda_m}$ is an orthonormal basis of $S_m$.
Notice that, if we set 
\begin{eqnarray*}
\phi_m^{(2)}&=&\frac{1}{{D_m}}\underset{t\in S_m^{(2)}\backslash\{0\}}
{\sup}\frac{\|t\|_\infty}{\|t\|},\\
\end{eqnarray*}
hypothesis M2 implies that $\phi_m^{(2)}$ is bounded by $r_0^2$.
The condition M1 must be replaced by the following condition: 
\begin{enumerate}
\item[M1'.] Each $S_m^{(2)}$ is a linear subspace of $(L^\infty\cap L^2)([0,1]^2)$ with dimension 
\hbox{$D_m^2\leq\sqrt{n}$.}
\end{enumerate}
Let now $$\gamma_n^{(2)}(t)=\frac{1}{n-1}\sum_{i=1}^{n-1}\{\|t\|^2-2t(X_i,X_{i+1})\}.$$
We define as above 
$$\hat{g}_m=\underset{t\in S_m^{(2)}}{\arg\min} \gamma_{n}^{(2)}(t)$$
and $\hat{m}^{(2)}=\underset{m\in \mathcal{M}_n}{\arg\min} 
                    [\gamma_{n}^{(2)}(\hat{g}_m)+\pen^{(2)}(m)]$
where $\pen^{(2)}(m)$ is a penalty function which would be specified later.
Lastly, we set $\tilde{g}=\hat{g}_{\hat{m}^{(2)}}$. 

\begin{thm}\label{th2}
Let $X_n$ be a Markov chain which satisfies Assumptions A1-A2-A3'-A4-A5 and 
$(S_m)_{m\in\mathcal{M}_n}$ 
be a collection of models satisfying Assumptions M1'-M2-M3. Then the estimator defined by 
\begin{equation} 
\tilde{g}=\hat{g}_{\hat{m}^{(2)}} \quad\text{ where } \quad
\hat{m}^{(2)}=\underset{m\in \mathcal{M}_n}{\arg\min} [\gamma_{n}^{(2)}(\hat{g}_m)+\pen^{(2)}(m)],
\end{equation} with 
\begin{equation}\pen^{(2)}(m)=K^{(2)}\frac{D_m^2}{n} \text{ for some  } K^{(2)}> K_0^{(2)}\end{equation} 
(where $K_0^{(2)}$ is a constant depending on the chain)
satisfies
$$\mathbb{E}\|\tilde{g}-g\|^2 \leq 3\underset{m\in\mathcal{M}_n}{\inf}\{d^2(g,S_m^{(2)}) 
    +\pen^{(2)}(m)\}+\frac{C_1}{n}$$
where $C_1$ does not depend on $n$.
 \end{thm}
The constant $K_0^{(2)}$ in the penalty is similar to the constant $K_0$ in Theorem \ref{th1}
(replacing $r_0$ by $r_0^2$ and $\|f\|_\infty$ by $\|g\|_\infty$).
We refer the reader to Remark \ref{rempenalite} for considerations related to these constants.

\begin{cor}\label{coro3}
Let $X_n$ be a Markov chain which satisfies Assumptions A1-A2-A3'-A4-A5 and 
$(S_m)_{m\in\mathcal{M}_n}$ be a collection of models mentioned in
Section \ref{ex} (with the regularity of polynomials and wavelets larger than $\alpha-1$).
If $g$ belongs to $B_{2,\infty}^\alpha([0,1]^2)$, with $\alpha>1$, then 
$$\mathbb{E}\|\tilde{g}-g\|^2 =O(n^{-\frac{2\alpha}{2\alpha+2}})$$
\end{cor}
This rate of convergence is the minimax rate for density estimation in dimension 
2 in the case of i.i.d. random variables (see for instance \citet{ibragimovhasminski}). 
Let us now proceed to the estimation of the transition density.

\subsection {Estimation of $\pi$}
The estimator of $\pi$ is defined in the following way.
Let $$\tilde{\pi}(x,y)=\begin{cases}
\frac{\tilde{g}(x,y)}{\tilde{f}(x)}& \text{ if }|\tilde{g}(x,y)|\leq a_n |\tilde{f}(x)|\\
0&\text{ else }\end{cases}$$
with $a_n=n^\beta$ and $\beta<1/8$.\\
We introduce a new assumption: 
\begin{enumerate}
\item[A6.] There exists a positive constant $\chi$ such that $\forall x \in [0,1], \quad f(x)\geq \chi$.\\
\end{enumerate}

\begin{thm}\label{th3}
Let $X_n$ be a Markov chain which satisfies Assumptions A1-A2-A3'-A4-A5-A6 and 
$(S_m)_{m\in\mathcal{M}_n}$ be a collection of models mentioned in
Section \ref{ex} (with the regularity of polynomials and wavelets larger than $\alpha-1$).
We suppose that the dimension $D_m$ of the models is such that
$$\forall m\in\mathcal{M}_n \quad \ln n\leq D_m \leq n^{1/4}.$$ 
If $f$ belongs to $B_{2,\infty}^\alpha([0,1])$, with $\alpha>1/2$, then for $n$ large enough
\begin{itemize}
\item there exists $C_1$ and $C_2$ such that
$$\E\|\pi-\tilde{\pi}\|^2\leq C_1\E\|g-\tilde{g}\|^2+C_2\E\|f-\tilde{f}\|^2+o(\frac{1}{n})$$
\item if furthermore $g$ belongs to $B_{2,\infty}^\beta([0,1]^2)$ (with $\beta>1$), then 
$$\displaystyle \E\|\pi-\tilde{\pi}\|^2=O(\sup(n^{-\frac{2\beta}{2\beta+2}},
  n^{-\frac{2\alpha}{2\alpha+1}}))$$
\end{itemize}
\end{thm}

\citet{clem} proved that $n^{-{2\beta}/{(2\beta+2)}}$ is the minimax rate for $f$ and $g$ 
of same regularity~ $\beta$. 
Notice that in this case the procedure is adaptive and there is no logarithmic loss in the 
estimation rate contrary to the result of \citet{clem}.

But it should be remembered that we consider only the restriction of $f$ or $\pi$ since 
the observations are in a compact set.
And the restriction of the stationary density
to $[0,1]$ may be less regular than the restriction of the transition density. 
The previous procedure has thus the disadvantage that the resulting rate 
does not depend only on the regularity of $\pi$ but also on the one of $f$.

However, if the chain lives on $[0,1]$ and if $g$ belongs to 
$B_{2,\infty}^\beta([0,1]^2)$ (that is to say that we consider the regularity
of $g$ on its whole support and not only on the compact of the observations) 
then equality $f(y)=\int g(x,y)dx $ yields that $f$ belongs to $B_{2,\infty}^\beta([0,1])$ 
and then $\E\|\pi-\tilde{\pi}\|^2=O(n^{-\frac{2\beta}{2\beta+2}})$. Moreover, 
if $\pi$ belongs to $B_{2,\infty}^\beta([0,1]^2)$, 
formula $f(y)=\int f(x)\pi(x,y)dx$ implies that $f$ belongs to 
$B_{2,\infty}^\beta([0,1])$. Then, by using properties of Besov spaces (see \citet{runst} 
p.192), $g=f\pi$ belongs to $B_{2,\infty}^\beta([0,1]^2)$. 
So in this case of a chain with compact support the minimax rate is
achieved as soon as $\pi$ belongs to $B_{2,\infty}^\beta([0,1]^2)$ with $\beta>1$.




\section{Simulations}

The computation of the previous estimator is very simple. We use the following 
procedure in 3 steps:
\begin{description}
\item First step: 
\begin{itemize}
    \item For each $m$, compute $\gamma_{n}(\hat{f}_m) +\pen(m)$. Notice that
$\gamma_{n}(\hat{f}_m)=-\sum_{\lambda\in\Lambda_m}\hat{\beta}_\lambda^2$
where $\hat{\beta}_\lambda$ is defined by \eqref{beta} and is quickly computed.
    \item Select the argmin $\hat{m}$ of $\gamma_{n}(\hat{f}_m) +\pen(m)$.
    \item Choose $\tilde{f}=\sum_{\lambda\in\Lambda_{\hat{m}}}\hat{\beta}_\lambda\varphi_\lambda$.
\end{itemize} 

\item Second step:
\begin{itemize}
    \item For each $m$ such that $D_m^2\leq \sqrt{n}$ compute 
$\gamma_{n}^{(2)}(\hat{g}_m)+\pen^{(2)}(m)$, with 
$\gamma_{n}^{(2)}(\hat{g}_m)=-\sum_{\lambda,\mu\in\Lambda_m}\hat{a}_{\lambda,\mu}^2$
where $\hat{a}_{\lambda,\mu}=({1}/{n})\sum_{i=1}^{n}\varphi_\lambda(X_i)\varphi_\mu(X_{i+1})$.
    \item Select the argmin $\hat{m}^{(2)}$ of $\gamma_{n}^{(2)}(\hat{g}_m)+\pen^{(2)}(m)$.
    \item Choose $\tilde{g}(x,y)=\sum_{\lambda,\mu\in\Lambda_{\hat{m}^{(2)}}}\hat{a}_{\lambda,\mu}\varphi_\lambda(x)
\varphi_\mu(y)$.
\end{itemize}

\item Third step: Compute
$\tilde{\pi}(x,y)={\tilde{g}(x,y)}/{\tilde{f}(x)}$ if $|\tilde{g}(x,y)|\leq n^{1/10} |\tilde{f}(x)|$ and 
$0$ otherwise.
\end{description}

The bases are here adjusted with an affin transform in order to be defined on
the estimation interval $[c,d]$ instead of
$[0,1]$. We consider 2 different bases (see Section \ref{ex}): 
trigonometric basis and histogram basis.

We found that a good choice for the penalty functions is 
$\pen(m)=5{D_m}/{n}$ and $\pen^{(2)}(m)= 0.02{D_m^2}/{n}.$

We consider several kinds of Markov chains : 
\begin{itemize}

\item An autoregressive process denoted by AR and defined by:
$$X_{n+1}=aX_n+b+\varepsilon_{n+1}$$
where the $\varepsilon_{n+1}$ are independent and identical distributed random variables,
with centered Gaussian distribution with variance $\sigma^2$. For this process, 
the stationary distribution is a Gaussian with mean $b/(1-a)$ and variance $\sigma^2/(1-a^2)$. 
By denoting by $\varphi(z)=1/(\sigma\sqrt{2\pi})\exp(-z^2/2\sigma^2)$ the Gaussian 
density, the transition density can be written 
$\pi(x,y)=\varphi(y-ax-b).$
We consider the following parameter values : \begin{itemize}
\item[(i)] $a=2/3$, $b=0$, $\sigma^2=5/9$, estimated on $[-2,2]^2$.
The stationary density of this chain is the standard Gaussian distribution.
\item[(ii)] $a=0.5$, $b=3$, $\sigma^2=1$, and then the process is estimated on $[4,8]^2$.
\end{itemize}

\item A radial Ornstein-Uhlenbeck process (in its discrete version).
For $j=1,\dots,\delta$, we define the processes: 
$\xi^j_{n+1}=a\xi^j_n+\beta\varepsilon_n^j$ where the $\varepsilon_{n}^j$ are i.i.d. standard Gaussian.
The chain is then defined by $X_n=\sqrt{\sum_{i=1}^\delta(\xi_n^i)^2}.$
The transition density is given in \cite{chaleyatgenon06} where this process is studied in detail:
$$\pi(x,y)=\1_{y>0}\exp\left(-\frac{y^2+a^2x^2}{2\beta^2}\right)I_{\delta/2-1}\left(\frac{axy}{\beta^2}\right)
\frac{ax}{\beta^2}\left(\frac{y}{ax}\right)^{\delta/2}$$ 
and $I_{\delta/2-1}$ is the Bessel function with index $\delta/2-1$.
The invariant density is 
$f(x)=C\1_{x>0}\exp(-{x^2}/{2\rho^2})x^{\delta-1}$
with $\rho^2=\beta^2/(1-a^2)$ and $C$ such that $\int f=1$.
This process (with here $a=0.5$, $\beta=3$, $\delta=3$) is denoted by $\sqrt{\text{CIR}}$ since its square
is actually a Cox-Ingersoll-Ross process. The estimation domain for this process is $[2,10]^2$.

\item A Cox-Ingersoll-Ross process, which is exactly the square of the previous process.
It follows a Gamma density for invariant distribution 
with scale parameter $l=1/2\rho^2$ and shape parameter $a=\delta/2$. The transition density is
$$\pi(x,y)=\frac{1}{2\beta^2}\exp\left(-\frac{y+a^2x}{2\beta^2}\right)
I_{\delta/2-1}\left(\frac{a\sqrt{xy}}{\beta^2}\right)\left(\frac{y}{a^2x}\right)^{\delta/4-1/2}$$ 
The used parameters are the following:
\begin{itemize}
\item[(iii)] $a= 3/4$, $b=\sqrt{7/48} $ (so that $l=3/2$) and $\delta=4$, estimated on $[0.1,3]^2$.
\item[(iv)] $a=1/3$, $b=3/4$ and $\delta=2$. This chain is estimated on $[0,2]^2$.
\end{itemize}

\item  An ARCH process defined by $X_{n+1}=\sin(X_n)+(\cos(X_n)+3)\varepsilon_{n+1}$
where the $\varepsilon_{n+1}$ are i.i.d. standard Gaussian.
The transition density of this chain is 
$$\pi(x,y)=\varphi\left(\frac{y-\sin(x)}{\cos(x)+3}\right)\frac{1}{\cos(x)+3}$$
and we estimate this process on  $[-5,5]^2$.
\end{itemize} 

For this last chain, the stationary density is not explicit. So we simulate $n+500$ variables and 
we estimate only from the last $n$ to ensure the stationarity of the process. For the other chains,
it is sufficient to simulate an initial variable $X_0$ with density $f$.

\begin{figure}
\includegraphics[scale=0.55]{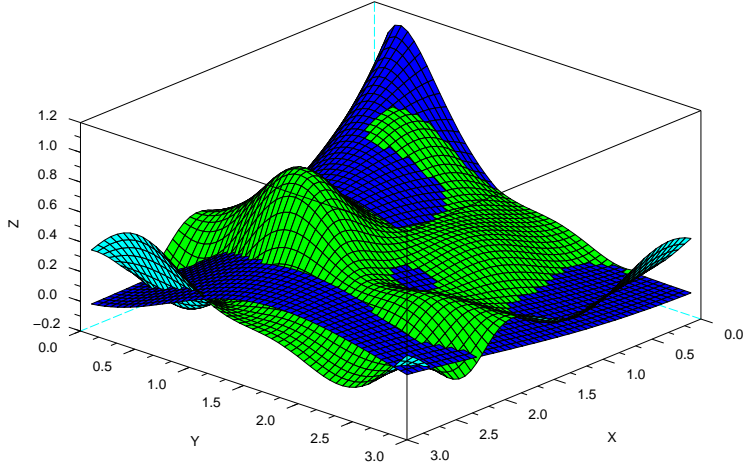}
\caption{Estimator (light surface) and true transition (dark surface) for the process CIR(iii) 
estimated with a trigonometric basis, n=1000}
\label{fig1}       
\end{figure}

Figure \ref{fig1} illustrates the performance of the method  and Table \ref{risqueL2p} shows the 
 $L^2$-risk for different values of $n$.

\begin{table}[!h]
\begin{tabular}{llllllc}
\hline\noalign{\smallskip}
{$n$} 
& 50 & 100 & 250 & 500 & 1000 & basis \\
\noalign{\smallskip}\hline
\hline
AR(i) & 0.7280 & 0.5442 & 0.2773 & 0.1868 & 0.1767 & H \\
& 0.5262 & 0.4682 & 0.2223 & 0.1797 & 0.1478 & T\\
\hline 
AR(ii) & 0.4798 & 0.3252 & 0.2249 & 0.1160 & 0.0842 & H \\
& 0.2867 & 0.2393 & 0.1770 & 0.1342 & 0.1083 & T\\
\hline 
$\sqrt{\text{CIR}}$  & 0.3054 & 0.2324 & 0.1724 & 0.1523 & 0.1278 & H \\
& 0.2157 & 0.1939 & 0.1450 & 0.1284 & 0.0815 & T\\
\hline 
CIR(iii) & 0.5086 & 0.3082 & 0.2113 & 0.1760 & 0.1477 & H\\
& 0.4170 & 0.3959 & 0.2843 & 0.2565 & 0.2265 & T\\
\hline 
CIR(iv) & 0.3381 & 0.2101 & 0.1205 & 0.0756 & 0.0458 & H\\
& 0.2273 & 0.2212 & 0.1715 & 0.1338 & 0.1328 & T\\
\hline 
ARCH & 0.3170 & 0.3013 & 0.2420 & 0.2124 & 0.1610 & H\\
& 0.2553 & 0.2541 & 0.2075 & 0.1884 & 0.1689 & T\\
\noalign{\smallskip}\hline
\end{tabular}
\caption{MISE $\E\|\pi-\tilde{\pi}\|^2$ averaged over 
$N=200$ samples. H: histogram basis, T: trigonometric basis.}
\label{risqueL2p}
\end{table}

The results in Table \ref{risqueL2p} are roughly good and illustrate that 
we can not pretend that a basis among the others gives better 
results. We can then imagine a mixed strategy, i.e. a procedure which uses several kinds of bases 
and which can choose the best basis or, for instance, the best degree for a polynomial basis. 
These techniques are successfully used in regression frameworks by  
\citet{comterozenholc2002,comterozenholc2004}.

The results for the stationary density are given in Table \ref{risqueL2f2}.

\begin{table}[!h]
\begin{tabular}{llllllc}
\hline\noalign{\smallskip}
{$n$} 
& 50 & 100 & 250 & 500 & 1000 & basis \\
\noalign{\smallskip}\hline
\hline
AR(i) & 0.0658 & 0.0599 & 0.0329 & 0.0137 & 0.0122 & H \\
& 0.0569 & 0.0538 & 0.0246 & 0.0040 & 0.0026 & T\\
\hline 
AR(ii) & 0.0388 & 0.0354 & 0.0309 & 0.0147 & 0.0081 & H \\
& 0.0342 & 0.0342 & 0.0327 & 0.0195 & 0.0054 & T\\
\hline 
$\sqrt{\text{CIR}}$  & 0.0127 & 0.0115 & 0.0105 & 0.0102 & 0.0096 & H \\
& 0.0169 & 0.0169 & 0.0168 & 0.0166 & 0.0107 & T\\
\hline 
CIR(iii) & 0.0335 & 0.0268 & 0.0229 & 0.0222 & 0.0210 & H\\
& 0.0630 & 0.0385 & 0.0216 & 0.0211 & 0.0191 & T\\
\hline 
CIR(iv) & 0.0317 & 0.0249 & 0.0223  & 0.0185 & 0.0103 & H\\
& 0.0873 & 0.0734 & 0.0572 & 0.0522 & 0.0458 & T\\
\noalign{\smallskip}\hline
\end{tabular}
\caption{MISE $\E\|f-\tilde{f}\|^2$ averaged over 
$N=200$ samples. H: histogram basis, T: trigonometric basis.}
\label{risqueL2f2}
\end{table}

We can compare results of Table \ref{risqueL2f2} with those of \cite{dalelane} who gives results of 
simulations for i.i.d. random variables. For density estimation, she uses three types of kernel: 
Gauss kernel, sinc-kernel (where $\operatorname{sinc}(x)=\sin (x)/x$) and her Cross Validation 
optimal kernel (denoted by Dal). Table \ref{dalelane} gives 
her results for the Gaussian density and the Gamma distribution with the same parameters that we 
used (2 and 3/2). If we compare the results that she obtains with 
her optimal kernel and our results with the trigonometric basis, we observe 
that her risks are about 5 times less than ours. However this kernel is particularly effective and 
if we consider the classical kernels, we notice that the results are almost comparable, with a 
reasonable price for dependency.

\begin{table}
\begin{tabular}{llllllll}
\hline\noalign{\smallskip}
$n$ & 100 & 500 & 1000 & kernel \\
\noalign{\smallskip}\hline\noalign{\smallskip}
& 0.0065 &  0.0013 & 0.0008 & Dal \\
Gaussian & 0.0127 &  0.0028 & 0.0016 & Gauss \\
(=AR(i)) & 0.0114 &  0.0026 & 0.0010 & sinc \\
\hline
 & 0.0148 & 0.0052 & 0.0027 & Dal \\
Gamma & 0.0209 & 0.0061 & 0.0031  & Gauss\\
(=CIR(iii)) & 0.0403 & 0.0166 & 0.0037 & sinc \\
 \noalign{\smallskip}\hline
\end{tabular}
\caption{MISE obtained by \citet{dalelane} for i.i.d. data, averaged over $50$ samples}
\label{dalelane}
\end{table}

\section {Proofs }
\subsection{The Nummelin splitting technique} \label{nummelin}
This whole subsection is summarized from \citet{loscher} p.60--63 and is detailed for the sake 
of completeness. 

The interest of the Nummelin splitting technique is to create a two-dimensional chain (the "split 
chain"), which contains automatically an atom. Let us recall the definition of an atom.
Let $A$ be a set such that $\psi(A)>0$ where $\psi$ is an irreducibility 
measure. The set $A$ is called an atom for the chain $(X_n)$ with transition 
kernel $P$ if there exists a measure $\nu$ such that $P(x,B)=\nu(B)$, for all $x$ in $A$ 
and for all event $B$.

Let us now describe the splitting method. Let $E=[0,1]$ the state space and $\mathcal{E}$ the 
associated $\sigma$-field. Each point $x$ in $E$ is splitted in $x_0=(x,0)\in E_0=E\times\{0\}$ 
and $x_1=(x,1)\in E_1=E\times\{1\}$. Each set $A$ in $\mathcal{E}$ is splitted in $A_0=A\times\{0\}$ 
and $A_1=A\times\{1\}$. Thus, we have defined a new probability space $(E^*, \mathcal{E}^*)$
where $E^*:=E_0\cup E_1$ and $\mathcal{E}^*=\sigma (A_0,A_1 : A\in \mathcal{E}).$ Using $h$ defined 
in A4, a measure $\lambda$ on $(E,\mathcal{E})$ splits according to 
$$\begin{cases}
\lambda^*(A_1)&=\int \1_A(x)h(x)\lambda(dx)\\
\lambda^*(A_0)&=\int \1_A(x)(1-h)(x)\lambda(dx)\\
\end{cases}$$
Notice that $\lambda^*(A_0\cup A_1)=\lambda(A).$ Now the aim is to define a new transition probability
$P^*(.,.)$ on $(E^*, \mathcal{E}^*)$ to replace the transition kernel $P$ of $(X_n)$. Let 
$$P^*(x_i,.)=\begin{cases}
\cfrac{1}{1-h(x)}(P-h\otimes\nu)^*(x,.)& \text{ if } i=0  \text{ and } h(x)>1 \\
\nu^* &   \text{ else }
\end{cases}$$
where $\nu$ is the measure introduced in A4 and $h\otimes\nu$ is a kernel defined by 
$h\otimes\nu(x,dy)=h(x)\nu(dy)$.
Consider now a chain $(X_n^*)$ on $(E^*, \mathcal{E}^*)$ with one-step transition $P^*$ and with 
starting law $\mu^*$. The split chain $(X_n^*)$ has the following properties:
 
\begin{enumerate}
\item[P1.]For all $(A_p)_{0\leq p\leq N}\in \mathcal{E}^N$ and for all measure $\lambda$
$$P_\lambda(X_p\in A_p, 0\leq p\leq N )= P_{\lambda^*}(X_p^*\in A_p\times\{0,1\}, 0\leq p\leq N ). $$
\item[P2.]The split chain is irreducible positive recurrent with stationary distribution $\mu^*$.
\item[P3.]The set $E_1$ is an atom for $(X_n^*)$.
\end{enumerate}

We can also extend functions $g: E\mapsto \mathbb{R}$ to $E^*$ via
$g^*(x_0)=g(x)=g^*(x_1).$
Then, the  property P1 can be written: for all function $\mathcal{E}$-measurable 
$g: E^N \mapsto \mathbb{R}$ $$\E_\lambda(g(X_1,.., X_N))=\E_{\lambda^*}(g^*(X_1^*,.., X_N^*)).$$
We can say that $(X_n)$ is a marginal chain of $(X_n^*)$.
When necessary, the following proofs are decomposed in two steps: first, we assume that the Markov 
chain has an atom, next we extend the result to the general chain by introducing the artificial atom 
$E_1$.

\subsection{Proof of Proposition \ref{est1}}
\emph{First step: } We suppose that $(X_n)$ has an atom $A$.

Let $f_m$ be the orthogonal projection of $f$ on $S_m$. Pythagoras theorem gives us: 
$$\E\|f-\hat{f}_m\|^2=d^2(f,S_m)+\E\|f_m-\hat{f}_m\|^2.$$
We recognize in the right member a bias term and a variance term. 
According to the expresssion \eqref{beta} of $\hat{f}_m$ the variance term can be written:  
\begin{equation}
\E\|f_m-\hat{f}_m\|^2=\sum_{\lambda\in\Lambda_m}\text{Var}(\hat{\beta}_{\lambda})
  =\sum_{\lambda\in\Lambda_m}\E(\nu_n^2(\varphi_{\lambda}))\label{varianceterm}\end{equation}
where $\nu_n(t)=(1/n)\sum_{i=1}^{n}[t(X_i)-<t,f>]$.
By denoting $\tau=\tau(1)=\inf\{n\geq 1, X_n\in A\}$ and 
$\tau(j)=\inf\{n>\tau(j-1), X_n\in A\}$ for $j\geq 2$, we can decompose 
$\nu_n(t)$ in the classic following way: 

\begin{equation}
\nu_n(t)=\nu_n^{(1)}(t)+\nu_n^{(2)}(t)+\nu_n^{(3)}(t)+\nu_n^{(4)}(t)\label{decompo}
\end{equation}

$\begin{disarray}{rrcl}
\text{with  }&\nu_n^{(1)}(t)&=&\nu_n(t) \1_{\tau >n} , \\
&\nu_n^{(2)}(t)&=&\frac{1}{n}\sum_{i=1}^{\tau}[t(X_i)-<t,f>]
                 \1_{\tau\leq n} ,\\
&\nu_n^{(3)}(t)&=&\frac{1}{n}\sum_{i=1+\tau(1)}^{\tau(l_n)}[t(X_i)-<t,f>]
                 \1_{\tau\leq n} ,\\
&\nu_n^{(4)}(t)&=&\frac{1}{n}\sum_{i=\tau(l_n)+1}^{n}[t(X_i)-<t,f>]
                 \1_{\tau\leq n} ,\\
\end {disarray}$\\

and $l_n=\sum_{i=1}^{n}\1_A(X_i)$ (number of visits to the atom $A$).
Hence, \begin{equation*}
{\nu_n}(t)^2\leq 4\{{\nu_n}^{(1)}(t)^2+{\nu_n}^{(2)}(t)^2+
   {\nu_n}^{(3)}(t)^2+{\nu_n}^{(4)}(t)^2 \}.
\end{equation*}

\noindent$\bullet$
To bound $\nu_n^{(1)}(t)^2$, notice that $|\nu_n(t)|\leq 2\|t\|_\infty.$
And then, by using M2 and \eqref{M2},
$|\nu_n^{(1)}(t)|\leq 2r_0\sqrt{D_m}\|t\|\1_{\tau >n}.$
Thus,
\begin{eqnarray*}
\E(\nu_n^{(1)}(t)^2)&\leq& 4r_0^2\|t\|^2D_mP(\tau>n)\leq
4r_0^2\|t\|^2\E(\tau^2)\frac{D_m}{n^2}.
\end{eqnarray*}

 \noindent$\bullet$
We bound the second term in the same way. Since
$|\nu_n^{(2)}(t)|\leq 2({\tau}/{n)}\|t\|_\infty$, we obtain
$|\nu_n^{(2)}(t)|\leq2\|t\|r_0\tau\sqrt{D_m}/n$ and then 
$$\E(\nu_n^{(2)}(t)^2)\leq 4r_0^2\|t\|^2\E(\tau^2)\frac{D_m}{n^2}.$$

  \noindent$\bullet$
Let us study now the fourth term. As
\begin{equation*}
|\nu_n^{(4)}(t)|\leq 2\frac{n-\tau(l_n)}{n}\|t\|_\infty\1_{\tau\leq n}
\leq2(n-\tau(l_n))\frac{\sqrt{D_m}}{n}r_0\|t\|\1_{\tau\leq n} ,
\end {equation*} we get
$\quad\E(\nu_n^{(4)}(t)^2)\leq 4r_0^2\|t\|^2\cfrac{D_m}{n^2}\E((n-\tau(l_n))^2\1_{\tau\leq n}).$\\
It remains to bound $\E((n-\tau(l_n))^2\1_{\tau\leq n})$:
\begin{eqnarray*}
&&\E_\mu((n-\tau(l_n))^2\1_{\tau\leq n})=\sum_{k=1}^n\E_\mu((n-k)^2\1_{\tau(l_n)=k}
  \1_{\tau\leq n})\\
&&=\sum_{k=1}^n(n-k)^2P_\mu(X_{k+1}\notin A,..,X_n\notin A|X_k\in A)P_\mu(X_k \in A)\\
&&=\sum_{k=1}^n(n-k)^2P_A(X_1\notin A,..,X_{n-k}\notin A)\mu(A)\\
\end{eqnarray*}
 by using the stationarity of $X$ and the Markov property. Hence 
\begin{eqnarray}
\E_\mu((n-\tau(l_n))^2\1_{\tau\leq n})&=&\sum_{k=1}^n(n-k)^2P_A(\tau>n-k)\mu(A)\nonumber\\
&\leq&\sum_{k=1}^{n-1}\frac{\E_A(\tau^4)}{(n-k)^2}\mu(A).\nonumber
\end{eqnarray} Therefore 
$\E_\mu((n-\tau(l_n))^2\1_{\tau\leq n})\leq2\E_A(\tau^4)\mu(A).$
Finally
$$\mathbb{E}(\nu_n^{(4)}(t)^2)\leq 8r_0^2\|t\|^2\mu(A)\mathbb{E}_A(\tau^4)\frac{D_m}{n^2}$$
and we can summarize the last three results by
\begin{equation}
\E\left(\nu_n^{(1)}(t)^2+\nu_n^{(2)}(t)^2+\nu_n^{(4)}(t)^2\right)
\leq 8r_0^2\|t\|^2[\E_\mu(\tau^2)+\mu(A)\E_A(\tau^4)]\frac{D_m}{n^2}.\label{124}
\end{equation}
In particular, if $t=\varphi_\lambda$, using that $D_m\leq n$, 
$$\E\left(\nu_n^{(1)}(\varphi_\lambda)^2+\nu_n^{(2)}(\varphi_\lambda)^2+
\nu_n^{(4)}(\varphi_\lambda)^2\right)\leq 8r_0^2\frac{\E_\mu(\tau^2)+\mu(A)\E_A(\tau^4)}{n}.$$

 \noindent$\bullet$
Last we can write 
$\nu_n^{(3)}(t)=(1/n)\sum_{j=1}^{l_n-1}S_j(t) \1_{\tau\leq n}$
where \begin{eqnarray}
S_j(t)=\sum_{i=1+\tau(j)}^{\tau(j+1)}(t(X_i)-<t,f>).\label{sj}
\end{eqnarray}
We remark that, according to the Markov property, the $S_j(t)$ are independent identically 
distributed and centered. Thus,
$$\mathbb{E}(\nu_n^{(3)}(\varphi_{\lambda})^2)\leq\frac{1}{n^2}\sum_{j=1}^{l_n-1}
\E|S_j(\varphi_{\lambda})|^2.$$
Then, we use Lemma \ref{maj1} below to bound the expectation of
$\nu_n^{(3)}(\varphi_{\lambda})^2$ :

\begin{lem}\label{maj1} 
For all $ m\geq 2$,
$\E_\mu|S_j(t)|^m \leq (2\|t\|_\infty)^{m-2}\|f\|_\infty\|t\|^2\E_A(\tau^m).$
\end{lem}

We can then give the bound
$$\mathbb{E}(\nu_n^{(3)}(\varphi_{\lambda})^2)\leq\frac{1}{n^2}\sum_{j=1}^{n}
\|f\|_\infty\|\varphi_{\lambda}\|^2\E_A(\tau^2)
\leq\frac{\|f\|_\infty\E_A(\tau^2)}{n}.$$
Finally
$$\E(\nu_n^2(\varphi_{\lambda}))\leq\frac{4}{n}[8r_0^2(\E_\mu(\tau^2)+\mu(A)\E_A(\tau^4))
+\|f\|_\infty\E_A(\tau^2)].$$
Let $C=4[8r_0^2(\E_\mu(\tau^2)+\mu(A)\E_A(\tau^4))+\|f\|_\infty\E_A(\tau^2)]$. We obtain
with \eqref{varianceterm}
$$\E\|f_m-\hat{f}_m\|^2\leq C\frac{D_m}{n}.$$

\emph{Second step: } We do not suppose any more that $(X_n)$ has an atom.

Let us apply the Nummelin splitting technique to the chain $(X_n)$ and let 
\begin{equation}
\gamma_{n}^*(t)=\frac{1}{n}\sum_{i=1}^{n}[\|t\|^2-2t^*(X_i^*)].\label{gamman*}
\end{equation}
We define also \begin{equation}
\hat{f}_m^*=\underset{t\in S_m}{\arg\min} \gamma_{n}^*(t).\label{est*}\end{equation}
Then the property P1 in Section \ref{nummelin} yields
$ \E\|f-\hat{f}_m^*\|^2= \E\|f-\hat{f}_m\|^2.$
The split chain having an atom (property P3), we can use the first step to deduce
$\E\|f-\hat{f}_m^*\|^2\leq d^2(f,S_m)+C{D_m}/{n}.$
It follows that $$\E\|f-\hat{f}_m\|^2\leq d^2(f,S_m)+C{D_m}/{n}.$$
\findemo\\

\textbf{Proof of Lemma \ref{maj1}: }\label{preuvemaj1}
For all $j$, 
$\E_\mu|S_j(t)|^m = \E_\mu|S_1(t)|^m=\E_\mu|\sum_{i=\tau+1}^{\tau(2)}\bar{t}(X_i)|^m$
where $\bar{t} =t-<t,f>$. Thus
\begin{eqnarray*}
&&\E_\mu|S_j(t)|^m=\sum_{k<l}\E\bigg(\Big|\sum_{i=k+1}^{l}\bar{t}(X_i)\Big|^m|\tau=k, \tau(2)=l\bigg)
        P(\tau=k, \tau(2)=l)\\
&&\leq\sum_{k<l}(2\|t\|_\infty(l-k))^{m-2}\E\bigg(\Big|\sum_{i=k+1}^{l}\bar{t}(X_i)\Big|^2
         |\tau=k, \tau(2)=l\bigg) P(\tau=k, \tau(2)=l)\\
&&\leq\sum_{k<l}(2\|t\|_\infty)^{m-2}(l-k)^{m-1}\sum_{i=k+1}^{l}\E\bigg(\big|\bar{t}(X_i)\big|^2
         |\tau=k, \tau(2)=l\bigg) P(\tau=k, \tau(2)=l)\\
\end{eqnarray*}
using the Schwarz inequality. Then, since the $X_i$ have the same distribution under $\mu$. 
\begin{eqnarray*}
\E_\mu|S_j(t)|^m&\leq&\sum_{k<l}(2\|t\|_\infty)^{m-2}(l-k)^{m}\E({t}^2(X_1))
        P(\tau=k, \tau(2)=l)\\
&\leq&\sum_{k<l}(2\|t\|_\infty)^{m-2}(l-k)^{m}\|f\|_\infty\|t\|^2
        P(\tau=k, \tau(2)=l)\\
&\leq&(2\|t\|_\infty)^{m-2}\E(|\tau(2)-\tau|^{m})\|f\|_\infty\|t\|^2.
\end{eqnarray*}
We conclude by using the Markov property.
\findemo

\subsection{Proof of Corollary \ref{coro1}}
According to Proposition \ref{est1}
$\E\|f-\hat{f}_m\|^2\leq d^2(f,S_m)+C{D_m}/{n}.$
Then we use Lemma 12 in \citet{BBM} which ensures that (for piecewise polynomials or wavelets having
a regularity larger than $\alpha-1$ and for trigonometric polynomials) 
$d^2(f,S_m) =O(D_m^{-2\alpha}).$ Thus,
$$\E\|f-\hat{f}_m\|^2=O(D_m^{-2\alpha}+ \frac{D_m}{n})$$
In particular, if  $D_m=\lfloor n^\frac{1}{1+2\alpha}\rfloor$,
then $ \E\|f-\hat{f}_m\|^2=O(n^{-\frac{2\alpha}{1+2\alpha}}).$
\findemo

\subsection{Proof of Theorem \ref{th1}}
\emph{First step: } We suppose that $(X_n)$ has an atom $A$.

Let $m$ in $\mathcal{M}_n $.  The definition of $\hat{m}$ yields that
$\gamma_n(\hat{f}_{\hat{m}})+\pen(\hat{m})\leq\gamma_n(f_m)+\pen(m) .$
This leads to \begin{equation}
\|\hat{f}_{\hat{m}}-f\|^2 \leq \|f_{m}-f\|^2 + 2\nu_n(\hat{f}_{\hat{m}}-f_m)
            +\pen(m)-\pen(\hat{m})\label{AA}
\end{equation}
where $\nu_n(t)=(1/n)\sum_{i=1}^{n}[t(X_i)-<t,f>]$.

\begin{rem}
If $t$ is deterministic, $\nu_n(t)$ can actually be written 
$\nu_n(t)=(1/n)\sum_{i=1}^{n}[t(X_i)-\E(t(X_i))]$.
\end{rem}

We set $B(m,m')=\{t\in S_m+S_{m'}, \|t\|=1\}$.
Let us write now
$$\begin{disarray}{c}
2\nu_n(\hat{f}_{\hat{m}}-f_m)=2\|\hat{f}_{\hat{m}}-f_m\|
  \nu_n\big(\frac{\hat{f}_{\hat{m}}-f_m}{\|\hat{f}_{\hat{m}}-f_m\|}\big)\\
\leq2\|\hat{f}_{\hat{m}}-f_m\|\underset{t\in B(m,\hat{m})}{~\sup~}\nu_n(t)
\leq\frac{1}{5}\|\hat{f}_{\hat{m}}-f_m\|^2+5\underset{t\in B(m,\hat{m})}{~\sup~}\nu_n(t)^2
\end{disarray}$$
by using inequality $2xy\leq \cfrac{1}{5}x^2+5y^2$. Thus,
\begin{eqnarray}
2\E|\nu_n(\hat{f}_{\hat{m}}-f_m)|&\leq&\frac{1}{5}\E\|\hat{f}_{\hat{m}}-f_m\|^2 
 +5\E(\underset{t\in B(m,\hat{m})}{~\sup~}\nu_n(t)^2).\label{BB}
\end {eqnarray}
Consider decomposition \eqref{decompo} of $\nu_n(t)$ again and let 
\begin{equation}
  \label{Zn}
  Z_n(t)=\frac{1}{n}\sum_{j=1+\tau(1)}^{\tau(l_n)}[t(X_i)-<t,f>].
\end{equation}
 Since $|\nu_n^{(3)}(t)|\leq |Z_n(t)|$, we can write
\begin{eqnarray*}
\underset{t\in B(m,\hat{m})}{~\sup~}\nu_n^{(3)}(t)^2
&\leq&p(m,\hat{m})+\sum_{m'\in \mathcal{M}_n}
      [\underset{t\in B(m,m')}{\sup}Z_n(t)^2-p(m,m')]_+
\end {eqnarray*}
where $p(.,.)$ is a function 
specified in Proposition \ref{pr} on page \pageref{pr}. 
Then, the bound \eqref{124} combined with M1, \eqref{AA} and \eqref{BB} gives
\begin{eqnarray*}
\mathbb{E}\|\hat{f}_{\hat{m}}-f\|^2& \leq &\|f_{m}-f\|^2
  +\frac{1}{5}\mathbb{E}\|\hat{f}_{\hat{m}}-f_m\|^2 
  +160r_0^2\frac{\E(\tau^2)+\mu(A)\E_A(\tau^4)}{n}\\&&
  +20\sum_{m'\in \mathcal{M}_n}\mathbb{E}[\underset{t\in B(m,m')}{~\sup~}Z_n(t)^2-p(m,m')]_+\\&&
  +\mathbb{E}(20p(m,\hat{m})+\pen(m)-\pen(\hat{m})).
\end{eqnarray*}

We choose $\pen(m)$ such that $20p(m,m')\leq \pen(m)+\pen(m')$. Thus
$ 20p(m,\hat{m})+\pen(m)-\pen(\hat{m})\leq 2\pen(m)$.
Let \begin{equation}\label{Wm}
\quad W(m,m')=[\underset{t\in B(m,m')}{~\sup~}Z_n^2(t)-p(m,m')]_{+} .
\end{equation}
We use now the inequality $\cfrac{1}{5}(x+y)^2\leq \cfrac{1}{3}x^2+\cfrac{1}{2}y^2$
to deduce
$$\mathbb{E}\|\hat{f}_{\hat{m}}-f\|^2 \leq \frac{1}{3}\mathbb{E}\|\hat{f}_{\hat{m}}-f\|^2 
  +\frac{3}{2}\|f_m-f\|^2+20\sum_{m'\in \mathcal{M}_n}\mathbb{E}W(m,m')
  +2\pen(m)+\frac{C}{n}$$
and thus
$$ \mathbb{E}\|\hat{f}_{\hat{m}}-f\|^2 \leq \frac{9}{4}\|f_m-f\|^2 
    +30\sum_{m'\in \mathcal{M}_n}\mathbb{E}W(m,m')+3\pen(m)+\frac{3C}{2n}.$$

We need now to bound $\mathbb{E}W(m,m')$ to complete the proof. 
Proposition \ref{pr} below implies 
$$\E W(m,m')\leq K'e^{-D_{m'}}(r_0\vee 1)^2K_3\frac{1+ K_2\|f\|_\infty}{n}$$
where $K'$ is a numerical constant and $K_2, K_3$ depend on the chain and with \begin{eqnarray}
p(m,m')&=&K\frac{{\rm dim}(S_m+S_{m'})}{n}(r_0\vee 1)^2K_3(1+K_2\|f\|_\infty)\label{pm}.
\end{eqnarray}
The notation $a\vee b$ means $\max(a,b).$

Assumption M3 yields $\sum_{m'\in \mathcal{M}_n}e^{-D_{m'}}\leq \sum_{k\geq 1} e^{-k}
=1/(e-1)$. Thus, by summation on $m'$ in $\mathcal{M}_n$ 
$$\sum_{m'\in \mathcal{M}_n}\mathbb{E}W(m,m')\leq 
K'\frac{1}{e-1}(r_0\vee 1)^2K_3\frac{1+K_2\|f\|_\infty}{n}.$$

It remains to specify the penalty, which has to satisfy 
$20p(m,m')\leq \pen(m)+\pen(m').$
The value of $p(m,m')$ is given by (\ref{pm}), so we set
$$\pen(m)\geq 20K\frac{D_m}{n}(r_0\vee 1)^2K_3(1+K_2\|f\|_\infty ) $$

Finally
$$\forall m \hspace{1cm} \mathbb{E}\|\hat{f}_{\hat{m}}-f\|^2 \leq 3\|f_m-f\|^2 
    +3\pen(m)+\frac{C_1}{n}$$
where $C_1$ depends on $r_0,\|f\|_\infty,\mu(A),\E_\mu(\tau^2),\E_A(\tau^4),K_2, K_3 $.
Since it is true for all $m$, we obtain the result.

\emph{Second step: } We do not suppose any more that $(X_n)$ has an atom.

The Nummelin splitting technique allows us to create the chain $(X_n^*)$ and to define
$\gamma_{n}^*(t)$ and $\hat{f}_m^*$ as above by \eqref{gamman*},\eqref{est*}.
 Set now 
$$\hat{m}^*=\underset{m\in \mathcal{M}_n}{\arg\min} [\gamma_{n}^*(\hat{f}_m^*)
          +\pen(m)] $$
and $\tilde{f}^*=\hat{f}_{\hat{m}^*}^*.$
The property P1 in Section \ref{nummelin} gives
$ \E\|f-\tilde{f}\|^2= \E\|f-\tilde{f}^*\|^2.$
The split chain having an atom, we can use the first step to deduce
$$\E\|f-\tilde{f}^*\|^2\leq 3\underset{m\in\mathcal{M}_n}{\inf}\{d^2(f,S_m) 
    +\pen(m)\}+\frac{C_1}{n}.$$
And then the result is valid when replacing $\tilde{f}^*$ by $\tilde{f}$.

\findemo

\begin{prop} \label{pr}
Let $(X_n)$  be a Markov chain which satisfies A1--A5 and $(S_m)_{m\in\mathcal{M}_n}$ 
be a collection of models satisfying M1--M3. We suppose that $(X_n)$ has an atom $A$.
Let $Z_n(t)$ and $W(m,m')$ defined by \eqref{Zn} and \eqref{Wm} with
$$p(m,m')=K\frac{{\rm dim}(S_m+S_{m'})}{n}(r_0\vee 1)^2\frac{1+\|f\|_\infty\E_A(s^\tau)}{(\ln s)^2}$$
(where $K$ is a numerical constant and $s$ is a real depending on the chain). Then
$$\E W(m,m')\leq K'e^{-D_{m'}}(r_0\vee 1)^2\frac{1+\|f\|_\infty\E_A(s^\tau)}{(\ln s)^2n}$$
\end{prop}

\noindent\textbf{Proof of Proposition \ref{pr}: }
We can write
$Z_n(t)=(1/n)\sum_{j=1}^{l_n-1}S_j(t)$
where $S_j(t)$ is defined by \eqref{sj}.
According to Lemma \ref{maj1}: 
$\E_\mu|S_j(t)|^m \leq (2\|t\|_\infty)^{m-2}$ $\|f\|_\infty\|t\|^2\E_A(\tau^m).$
Now, we use condition A5 of geometric ergodicity. The proof of Theorem 15.4.2 in \citet{m&t} 
shows that $A$ is a Kendall set, i.e. there exists $s>1$ (depending on $A$) such that 
$\sup_{x\in A} \E_x(s^{\tau})<\infty$.  
Then $\E_A(\tau^m)\leq [m!/(\ln s)^m]\E_A(s^\tau)$. Indeed 
\begin{eqnarray*}
\E_A(\tau^m)&=&\int_0^\infty mx^{m-1}P_A(\tau>x)dx\\
&\leq&\int_0^\infty mx^{m-1}s^{-x}\E_A(s^\tau)dx
= \frac{m!}{(\ln s)^m}\E_A(s^\tau)
\end{eqnarray*}
Thus \begin{equation}
\forall m\geq 2 \hspace{1cm} \E_\mu|S_j(t)|^m \leq m!\left(\frac{2\|t\|_\infty}{\ln s}\right)^{m-2}
\frac{\|f\|_\infty\|t\|^2}{(\ln s)^2} \E_A(s^\tau).\label{maj2}
\end{equation}
We use now the following inequality (see \citet{petrov} p.49):
$$P\big(\underset{1\leq l\leq n}{\max} \sum_{j=1}^{l} S_j(t)\geq y\big)\leq 2P\big(\sum_{j=1}^{n} 
    S_j(t)\geq y-\sqrt{2B_n}\big)$$
where $B_n\geq \sum_{j=1}^{n}\E S_j(t)^2$.
The inequality (\ref{maj2}) gives us $B_n=2n\cfrac{\|f\|_\infty\|t\|^2}{(\ln s)^2} \E_A(s^\tau)$ and 
$$P\big(\sum_{j=1}^{l_n-1} S_j(t)\geq y\big)\leq P\big(\underset{1\leq l\leq n}{\max} 
\sum_{j=1}^{l} S_j(t)\geq y\big)\leq 2P\big(\sum_{j=1}^{n}S_j(t)\geq y-2\sqrt{n}\|t\|M/\ln s\big)$$
where $M^2=\|f\|_\infty \E_A(s^\tau)$.
We use then the Bernstein inequality given by \citet{birge&massart98}. 
$$P(\sum_{j=1}^{n}S_j(t)\geq n\varepsilon)\leq e^{-nx}$$
with $\varepsilon=\cfrac{2\|t\|_\infty}{\ln s}x+\cfrac{2\|t\|M}{\ln s}\sqrt{x}$ .
Indeed, according to (\ref{maj2}),  $$\frac{1}{n}\sum_{j=1}^{n}\E|S_j(t)|^m
  \leq \frac{m!}{2}(\frac{2\|t\|_\infty}{\ln s})^{m-2}(\frac{\sqrt{2}\|t\|M}{\ln s})^2.$$

Finally \begin{equation}\label{ineg}
P\left(Z_n(t)\geq \cfrac{2}{\ln s}\left[\|t\|_\infty x+M\|t\|\sqrt{x}
+M\|t\|/\sqrt{n}\right]\right)\leq 2e^{-nx}. 
\end{equation}

We will now use a chaining technique used in \citet{BBM}. Let us recall first the following lemma
(Lemma 9 p.400 in \citet{BBM}, see also Proposition 1 in \citet{birge&massart98}).

\begin{lem} \label{chainage}
Let $\bar S$ a subspace of $L^2$  with dimension $D$ spanned by $(\varphi_{\lambda})_
{\lambda\in\Lambda}$ (orthonormal basis). Let $$r=\frac{1}{\sqrt{D}}
\underset{\beta\neq 0}{\sup}\frac{\|\sum_{\lambda\in\Lambda}
 \beta_{\lambda}\varphi_{\lambda}\|_\infty}{\sup_{\lambda\in\Lambda}|\beta_{\lambda}|}.$$
Then, for all $\delta>0$, we can find a countable set $T\subset\bar S$ and a mapping 
$\pi$ from $\bar S$ to $T$ such that : 
\begin{itemize}
    \item for all ball $\mathcal{B}$ with radius $\sigma\geq 5\delta$
\begin{equation} |T\cap \mathcal{B}|\leq (5\sigma/\delta)^D    \label{9-1}
\end{equation}
    \item $\|u-\pi(u)\|\leq \delta$, $\forall u\in \bar S$ and 
$\sup_{u\in \pi^{-1}(t)}\|u-t\|_{\infty} \leq r\delta, ~~\forall t \in T .$
\end{itemize}
\end{lem}

 We apply this lemma to the subspace $S_m+ S_{m'}$ with dimension $D_m\vee D_{m'}$ denoted by 
$D(m,m')$ and $ r= r(m,m')$ defined by
$$r(m,m')=\frac{1}{\sqrt{D(m,m')}}\underset{\beta\neq 0}{\sup}\frac{\|\sum_{\lambda\in\Lambda(m,m')}
       \beta_{\lambda}\varphi_{\lambda}\|_\infty}{\sup_{\lambda\in\Lambda(m,m')}|\beta_{\lambda}|} $$
where $(\varphi_{\lambda})_{\lambda\in\Lambda(m,m')}$ is an orthonormal basis of $S_m+ S_{m'}$.
Notice that this quantity satisfy
$\phi_{m"}\leq r(m,m')\leq\sqrt{D(m,m')}\phi_{m"}$
where $m"$ is such that $S_m+S_{m'}=S_{m"}$ and then, using M2,
$$r(m,m')\leq r_0\sqrt{D(m,m')}.$$
We consider $\delta_0\leq 1/5$ , $\delta_k=\delta_02^{-k}$,
and the $T_k=T\cap B(m,m')$ where $T$ is defined by Lemma \ref{chainage} with $\delta=\delta_k$
and $ B(m,m')$ is the unit ball of $ S_m+S_{m'}$.
Inequality (\ref{9-1}) gives us $|T\cap  B(m,m')|\leq (5/\delta_k)^{D(m,m')} $. 
By letting $H_k=\ln(|T_k|)$, we obtain 
\begin{equation}H_k\leq D(m,m')[\ln(\frac{5}{\delta_0})+k\ln2].\label{Hk}\end{equation}
Thus, for all $u$ in $B(m,m')$, we can find a sequence $\{u_k\}_{k\geq 0}$ with $u_k\in T_k$
such that $\|u-u_k\|\leq \delta_k$ and $\|u-u_k\|_{\infty}\leq  r(m,m') \delta_k$.
Hence, we have the following decomposition:
$$u=u_0+\sum_{k=1}^{\infty} (u_{k}-u_{k-1})$$
with $\|u_0\|\leq 1$ and $\|u_0\|_{\infty}\leq r_0\sqrt{D(m,m')}\|u_0\|\leq r_0\sqrt{D(m,m')}$
and for all $k\geq 1$,
\begin{eqnarray*}
\|u_{k}-u_{k-1}\|&\leq &\delta_k+ \delta_{k-1}=3\delta_{k-1}/2, \\
\|u_k-u_{k-1}\|_{\infty}&\leq& 3 r(m,m')\delta_{k-1}/2\leq 3 r_0\sqrt{D(m,m')}\delta_{k-1}/2.
\end{eqnarray*}
Then  
\begin{align*}
P(\underset{u\in B(m,m')}{\sup} Z_n(u)>\eta)=&P(\exists (u_k)_{k\geq 0}\in \prod_{k\geq 0}T_k , 
Z_n(u_0)+  \sum_{k=1}^\infty Z_n(u_k-u_{k-1})>\eta_0+\sum_{k=1}^\infty \eta_k)\\
\leq& \sum_{u_0\in T_0}P(Z_n(u_0)>\eta_0)+\sum_{k=1}^\infty\sum_{\substack{u_k\in T_k\\ 
 u_{k-1}\in T_{k-1}}}P(Z_n( u_k-u_{k-1})>\eta_k)
\end{align*}
with $\eta_0+\sum_{k=1}^\infty \eta_k\leq\eta$.
We use the exponential inequality (\ref{ineg}) to obtain
\begin{align*}
\sum_{u_0\in T_0}P(Z_n(u_0)>\eta_0)\leq &2e^{H_0-nx_0}\\
\sum_{\substack{u_k\in T_k\\u_{k-1}\in T_{k-1}}}P(Z_n( u_k-u_{k-1})>\eta_k)\leq &2e^{H_k+H_{k-1}-nx_k}
\end{align*} 
by choosing 
$\left\{\begin{array}{l}
\eta_0=\cfrac{2}{\ln s}\left(r_0\sqrt{{D(m,m')}}x_0 +M\sqrt{x_0}+\cfrac{M}{\sqrt{n}}\right) \\
\eta_k=\cfrac{3}{\ln s}\left(r_0\sqrt{D(m,m')} \delta_{k-1}x_k +M\delta_{k-1}\sqrt{x_k}+
       \cfrac{M\delta_{k-1}}{\sqrt{n}}\right) .
\end{array}\right. $\\
Let us choose now the $(x_k)_{k\geq 0}$ such that $nx_0=H_0+D_{m'}+v$ 
and for $k\geq 1$,$$nx_k=H_{k-1}+H_k+kD_{m'}+D_{m'}+v$$
Thus \begin{eqnarray*}
P(\underset{u\in B(m,m')}{\sup} Z_n(u)>\eta)&\leq& 2e^{-D_{m'}-v}(1+\sum_{k\geq 1}e^{-kD_{m'}})
\leq 3.2e^{-D_{m'}-v}
\end{eqnarray*}
It remains to bound $\sum_{k=0}^\infty\eta_k$:
\begin{eqnarray*}
\sum_{k=0}^\infty\eta_k&\leq &\frac{1}{(\ln s)}(A_1+A_2+A_3).
\end{eqnarray*}
where  $\left\{\begin{disarray}{rcl}
A_1&=&r_0\sqrt{{D(m,m')}}(2x_0+3\sum_{k=1}^\infty \delta_{k-1}x_k)\\
A_2&=&2M\sqrt{x_0}+3M\sum_{k=1}^\infty\delta_{k-1}\sqrt{x_k}\\
A_3&=&2\cfrac{M}{\sqrt{n}} +\sum_{k=1}^\infty\cfrac{3M\delta_{k-1}}{\sqrt{n}}\\
\end{disarray}\right.$

\noindent$\bullet$
Regarding the third term, just write
\begin{eqnarray*}
A_3&=&\cfrac{M}{\sqrt{n}}\Big(2+3\sum_{k=1}^\infty\delta_{k-1}\Big)
    =\cfrac{M}{\sqrt{n}}(6\delta_0+2)\leq c_1(\delta_0)\cfrac{M}{\sqrt{n}}
\end{eqnarray*}
with $\begin{disarray}{rcl}
c_1(\delta_0)&=&6\delta_0+2.\\
\end{disarray}$

\noindent$\bullet$
Let us bound the first term.
First, recall that $D(m,m')\leq \sqrt{n}$ and then 
\begin{eqnarray*}
A_1&\leq&r_0\sqrt{\frac{n}{D(m,m')}}\left(2\cfrac{H_0+D_{m'}+v}{n} +3\sum_{k=1}^\infty \delta_{k-1}
     \cfrac{H_{k-1}+H_k+kD_{m'}+D_{m'}+v}{n}\right).\\
\end{eqnarray*}
Observing that $\sum_{k=1}^\infty \delta_{k-1}=2\delta_0$ and 
$\sum_{k=1}^\infty k\delta_{k-1}=4\delta_0$ and using \eqref{Hk}, we get
\begin{eqnarray*}
A_1&\leq&c_1(\delta_0)r_0 \frac{v}{\sqrt{nD(m,m')}} +
    c_2(\delta_0)r_0\sqrt{\frac {D(m,m')}{n}}\end{eqnarray*}
with 
$c_2(\delta_0)=c_1(\delta_0)+\ln(5/\delta_0)(2+12\delta_0)
     +6\delta_0(2+3\ln2)$

\noindent$\bullet$
To bound the second term, we use the Schwarz inequality and the inequality 
$\sqrt{a+b}\leq\sqrt{a}+\sqrt{b}$. We obtain
\begin{eqnarray*}
A_2&\leq&c_1(\delta_0)M\sqrt{\frac{v}{n}}+c_3(\delta_0)M\sqrt{\frac{D(m,m')}{n}}\\
\end{eqnarray*}
with 
$c_3(\delta_0)=2\sqrt{1+\ln(5/\delta_0)}+3\sqrt{2\delta_0}\sqrt{
                6\delta_0(1+\ln 2)+4\delta_0\ln(5/\delta_0)}$

We get so \begin{eqnarray*}
(\sum_{k=0}^\infty \eta_k)&\leq&\left(\frac{r_0\vee 1}{\ln s}\right)c_1
  \left(\frac{v}{\sqrt{nD(m,m')}}+M\sqrt{\frac{v}{n}}\right)\\ &&
  +\sqrt{\frac{D(m,m')}{n}}\left(\frac{r_0\vee 1}{\ln s}\right)
  [c_2+c_3M+c_1M]\\
(\sum_{k=0}^\infty \eta_k)^2&\leq& c_4(\delta_0)\left(\frac{r_0\vee 1}{\ln s}\right)^2
   [\frac{v^2}{nD(m,m')}\vee M^2\frac{v}{n}]\\ &&
   +c_5(\delta_0)\frac{D(m,m')}{n}\left(\frac{r_0\vee 1}{\ln s}\right)^2(1+M)^2
\end{eqnarray*}
where $\begin{cases}
c_4(\delta_0)=&6c_1^2\\
c_5(\delta_0)=&({6}/{5})\sup(c_2,c_3+c_1)^2
\end{cases}$\\
Let us choose now $\delta_0=0.024$ and then $c_4=28$, $c_5=268$.
Let $K_1=c_4(r_0\vee 1/\ln s)^2 $. Then
$$\eta^2=K_1[\frac{v^2}{nD(m,m')}\vee M^2\frac{v}{n}]+p(m,m')$$
where
$$p(m,m')=c_5(r_0\vee 1)^2\frac{D(m,m')}{n}\frac{1+\|f\|_\infty\E_A(s^\tau)}{(\ln s )^2}$$

$\begin{disarray}{llll}
\text{We get  }
&&P(\underset{u\in B(m,m')}{\sup} Z_n^2(u)>K_1 [\frac{v^2}{nD(m,m')}\vee M^2\frac{v}{n}]+p(m,m'))\\
&=& P(\underset{u\in B(m,m')}{\sup} Z_n^2(u)>\eta^2)\\
&\leq&P(\underset{u\in B(m,m')}{\sup} Z_n(u)>\eta)+P(\underset{u\in B(m,m')}{\sup} Z_n(u)<-\eta)\\
\end{disarray}$\\
Now 
\begin{align*}
P(\underset{u\in B(m,m')}{\sup} Z_n(u)<-\eta)\leq& \sum_{u_0\in T_0}P(Z_n(u_0)<-\eta_0)
  +\sum_{k=1}^\infty\sum_{\substack{u_k\in T_k\\u_{k-1}\in T_{k-1}}}P(Z_n( u_k-u_{k-1})<-\eta_k)\\
\leq& \sum_{u_0\in T_0}P(Z_n(-u_0)>\eta_0)
  +\sum_{k=1}^\infty\sum_{\substack{u_k\in T_k\\u_{k-1}\in T_{k-1}}}P(Z_n(-u_k+u_{k-1})>\eta_k)\\
\leq& 3.2e^{-D_{m'}-v}.
\end{align*}
 
Hence $$P(\underset{u\in B(m,m')}{\sup} Z_n^2(u)>K_1
     [\frac{v^2}{nD(m,m')}\vee M^2\frac{v}{n}]+p(m,m'))\leq 6.4e^{-D_{m'}-v}.$$

We obtain then 
$$\begin{disarray}{l}
\E[\underset{t\in B(m,m')}{~\sup~}Z_n^2(t)-p(m,m')]_{+}
\leq\int_0^\infty P(\underset{u\in B(m,m')}{\sup} Z_n^2(u)>p(m,m')+z)dz\\
\leq \int_0^{M^2D(m,m')}P(\underset{u\in B(m,m')}{\sup} Z_n^2(u)>p(m,m')+
    K_1M^2\frac{v}{n})K_1\frac{M^2}{n}dv\\
    +\int_{M^2D(m,m')}^\infty P(\underset{u\in B(m,m')}{\sup} Z_n^2(u)>p(m,m')
    +K_1\frac{v^2}{nD(m,m')})K_1\frac{2v}{nD(m,m')}dv\\
\leq\frac{K_1}{n}\left[M^2\int_0^\infty6.4e^{-D_{m'}-v}dv
    +\frac{2}{D(m,m')}\int_0^\infty6.4e^{-D_{m'}-v}vdv\right]\\
\leq\frac{6.4K_1}{n}e^{-D_{m'}}(M^2+\frac{2}{D(m,m')})
\leq12.8K_1e^{-D_{m'}}\frac{1+M^2}{n}.
\end{disarray}$$
By replacing $M^2$ by its value, we get so
$$\E W(m,m')\leq K'(\frac{r_0\vee 1}{\ln s})^2e^{-D_{m'}}\frac{1+\|f\|_\infty\E_A(s^\tau)}{n}$$
where $K'$ is a numerical constant
\findemo

\subsection{Proof of Corollary \ref{coro2}}
According to Theorem \ref{th1},
$\mathbb{E}\|\tilde{f}-f\|^2 \leq C_2\underset{m\in\mathcal{M}_n}{\inf}\{d^2(f,S_m) 
    +{D_m}/{n}\}.$
Since $d^2(f,S_m) =O(D_m^{-2\alpha})$ (see Lemma 12 in \citet{BBM}), 
$$\mathbb{E}\|\tilde{f}-f\|^2 \leq C_3\underset{m\in\mathcal{M}_n}{\inf}\{D_m^{-2\alpha} 
    +\frac{D_m}{n}\}$$
In particular, if $m_0$ is such that $D_{m_0}=\lfloor n^\frac{1}{1+2\alpha}\rfloor$,
then $$\mathbb{E}\|\tilde{f}-f\|^2 \leq C_3\{D_{m_0}^{-2\alpha} +\frac{D_{m_0}}{n}\}
\leq C_4n^{-\frac{2\alpha}{1+2\alpha}}.$$
The condition $D_m\leq \sqrt{n}$ allows this choice of $m$ only if $\alpha>\frac{1}{2}$.
\findemo

\subsection{Proof of Theorem \ref{th2}}
The proof is identical to the one of Theorem \ref{th1}.
\findemo

\subsection{Proof of Corollary \ref{coro3}}
It is sufficient to prove that $d(g,S_m^{(2)})\leq D_m^{-\alpha}$ if $g$ belongs to 
$B_{2,\infty}^\alpha([0,1]^2)$. It is done in the following lemma.
\findemo

\begin{lem}\label{approx}
Let $g$ in the Besov space $B_{2,\infty}^\alpha([0,1]^2)$. We consider 
the following spaces of dimension $D^2$ :\begin{itemize}
\item $S_1$ is a space of piecewiwe polynomials of degree bounded by $s>\alpha -1$
based on a partition with square of vertice $1/D$,
\item $S_2$ is a space of of orthonormal wavelets of regularity $s>\alpha -1$, 
\item $S_3$ is the space of trigonometric polynomials. 
\end{itemize}
Then, there exist positive constants $C_i$ such that
$$d(g,S_i)\leq C_i D^{-\alpha}\quad\text{for } i=1,2,3 .$$ 
\end{lem}

\emph{Proof of Lemme \ref{approx}: }
Let us recall the definition of $B_{2,\infty}^\alpha([0,1]^2)$. Let 
$$\Delta_h^rg(x,y)=\sum_{k=0}^r(-1)^{r-k}\binom{r}{k}g(x+kh_1,y+kh_2)$$
the rth difference operateur with step $h$ and 
$\omega_r(g,t)=\underset{|h|\leq t}{\sup}\|\Delta_h^rg\|_2$
the rth modulus of smoothness of g.
We say $g$ is in the Besov space $B_{2,\infty}^\alpha([0,1]^2)$ if 
 $\sup_{t>0}t^{-\alpha}\omega_r(g,t)<\infty$
for $r=\lfloor\alpha\rfloor +1$, or equivalently, for $r$ an integer larger than $\alpha$.

\citet{nlacta} proved that $d(g,S_1)\leq C \omega_{s+1}(g,D^{-1})$ , so
$$d(g,S_1)\leq C D^{-\alpha}.$$

For the wavelets case, we use the fact that $f$ belongs to $B_{2,\infty}^\alpha([0,1]^2)$ 
if and only if $\underset{j\geq-1}{\sup}2^{j\alpha}\|\beta_j\|<\infty$
(see \citet{4} chapter 6, section 10).
If $g_D$ is the orthogonal projection of $g$ on $S_2$, it follows from Bernstein's 
inequality that
$$\|g-g_D\|^2=\sum_{j>m}\sum_{k,l}|\beta_{jkl}|^2
       \leq C\sum_{j>m}2^{-2j\alpha}\leq C'D^{-j\alpha}$$
where $m$ is such that $2^m=D$.

For the trigonometric case, it is proved in \cite{Nikolskij} (p. 191 and 200) that
$d(g,S_3)\leq C \omega_{s+1}(g,D^{-1})$ so that $d(g,S_3)\leq C'D^{-\alpha}.$
 \findemo

\subsection{Proof of Theorem \ref{th3}}
Let us prove first the first item.
Let $E_n=\{\|f-\tilde{f}\|_\infty\leq {\chi}/{2}\}$ and $E_n^c$ its complementary.
On $E_n$, $\tilde{f}(x)=\tilde{f}(x)-f(x)+f(x)\geq\chi/{2}$
and for $n$ large enough,
$\displaystyle\tilde{\pi}(x,y)=\frac{\tilde{g}(x,y)}{\tilde{f}(x)} .$
For all $(x,y)\in [0,1]^2$,
$$\begin{disarray}{rcl}
|\tilde{\pi}(x,y)-\pi(x,y)|^2&\leq &|\frac{\tilde{g}(x,y)-\tilde{f}(x)\pi(x,y)}
{\tilde{f}(x)}|^2\1_{E_n}
+(\|\tilde{\pi}\|_\infty+\|\pi\|_\infty)^2\1_{E_n^c}\\
&\leq &\frac{|\tilde{g}(x,y)-g(x,y)+\pi(x,y)(f(x)-\tilde{f}(x))|^2}
{\chi^2/4}\\&&+(a_n+\|\pi\|_\infty)^2\1_{E_n^c}\\
\E\|\pi-\tilde{\pi}\|^2&\leq &\frac{8}{\chi^2}[\E\|g-\tilde{g}\|^2
  +\|\pi\|_\infty^2\E\|f-\tilde{f}\|^2]+(a_n+\|\pi\|_\infty)^2P(E_n^c)
\end{disarray}$$
It remains to bound $P(E_n^c)$. To do this, we observe that
$$\|f-\tilde{f}\|_\infty \leq \|f-f_{\hat{m}}\|_\infty+\|f_{\hat{m}}-\hat{f}_{\hat{m}}\|_\infty$$
Let $\gamma=\alpha-\frac{1}{2}$, then $B_{2,\infty}^\alpha([0,1])\subset 
B_{\infty,\infty}^\gamma([0,1])$
(see \citet{devorelorentz} p.182). 
Thus $f$ belongs to $B_{\infty,\infty}^\gamma([0,1])$ and Lemma 12 in \citet{BBM} gives  
$$\|f-f_{\hat{m}}\|_\infty\leq D_{\hat{m}}^{-\gamma}\leq (\ln n)^{-\gamma}$$

Thus $\|f-f_{\hat{m}}\|_\infty$ decreases to $0$ and $\|f-f_{\hat{m}}\|_\infty\leq {\chi}/{4}$ 
for $n$ large enough. So
$$P(E_n^c)\leq P(\|f_{\hat{m}}-\hat{f}_{\hat{m}}\|_\infty>\frac{\chi}{4})$$
But $\|f_{\hat{m}}-\hat{f}_{\hat{m}}\|_\infty\leq r_0\sqrt{D_{\hat{m}}}
\|f_{\hat{m}}-\hat{f}_{\hat{m}}\|\leq r_0n^{1/8}\|f_{\hat{m}}-\hat{f}_{\hat{m}}\|$
and $\|f_{\hat{m}}-\hat{f}_{\hat{m}}\|^2=\sum_{\lambda\in\Lambda_{\hat{m}}}
\nu_n^2(\varphi_\lambda)$.
Thus,
$$\begin{disarray}{rcl}
P(E_n^c)&\leq& P(\sum_{\lambda\in\Lambda_{\hat{m}}}\nu_n^2(\varphi_\lambda)>
  \frac{\chi^2}{16r_0^2n^{1/4}})\\
&\leq& P(\sum_{\lambda\in\Lambda_{\hat{m}}}\nu_n^{(1)}(\varphi_\lambda)^2
  +\nu_n^{(2)}(\varphi_\lambda)^2+\nu_n^{(4)}(\varphi_\lambda)^2>\frac{\chi^2}{32r_0^2n^{1/4}})
  \\&&+P(\sum_{\lambda\in\Lambda_{\hat{m}}}
  Z_n^2(\varphi_\lambda)>\frac{\chi^2}{32r_0^2n^{1/4}})\\
&\leq& \frac{32r_0^2n^{1/4}}{\chi^2}\E(\sum_{\lambda\in\Lambda_{\hat{m}}}\nu_n^{(1)}
  (\varphi_\lambda)^2+\nu_n^{(2)}(\varphi_\lambda)^2+\nu_n^{(4)}(\varphi_\lambda)^2)
  \\&&+\underset{m\in\mathcal{M}_n}{\sup}\sum_{\lambda\in\Lambda_m}P(Z_n^2(\varphi_\lambda)>
   \frac{\chi^2}{32r_0^2n^{1/2}})\\
\end{disarray}$$
We need then to bound two terms.
For the first term, 
let $S_{m_0}$ the maximum model with cardinal $D_{m_0}\leq n^{1/4}$.
Since $\Lambda_{\hat{m}}\subset \Lambda_{m_0}$ and 
using inequality \eqref{124} and the assumption $\forall m \quad D_m\leq n^{1/4}$, we obtain 
$$\frac{32r_0^2n^{1/4}}{\chi^2}\E(\sum_{\lambda\in\Lambda_{\hat{m}}}\nu_n^{(1)}
  (\varphi_\lambda)^2+\nu_n^{(2)}(\varphi_\lambda)^2+\nu_n^{(4)}(\varphi_\lambda)^2)
  \leq C'n^{-5/4}$$
Besides, for all $x$ and for all $\lambda$, using \eqref{ineg},
$$P(Z_n(\varphi_\lambda)\geq 2r_0n^{1/8}x+2M\sqrt{x}+2\frac{M}{\sqrt{n}})
 \leq 2e^{-nx}$$
and so
$$P(Z_n^2(\varphi_\lambda)\geq (2r_0n^{1/8}x+2M\sqrt{x}+2\frac{M}{\sqrt{n}})^2)
 \leq 4e^{-nx}$$

Let now $x=n^{-3/4}$, $x$ verifies (for $n$ large enough)
$$2r_0n^{3/8}x+2Mn^{1/4}\sqrt{x}+2Mn^{-1/4}\leq \frac{\chi}{r_0\sqrt{32}}$$
that yields
$$(2r_0n^{1/8}x+2M\sqrt{x}+2\frac{M}{\sqrt{n}})^2\leq \frac{\chi^2}{32r_0^2n^{1/2}}$$
The previous inequality gives then 
$$P\left(Z_n^2(\varphi_\lambda\right)>\frac{\chi^2}{32r_0^2n^{1/4}})\leq 4e^{-nx}
\leq 4e^{-{n}^{1/4}}$$
Finally
$$P(E_n^c)\leq 4n^{1/4}e^{-{n}^{1/4}}+C'n^{-5/4}\leq C"n^{-5/4}$$ for $n$ great enough.
And then, for $n$ large enough, $(a_n+\|\pi\|_\infty)^2P(E_n^c)\leq Ca_n^2n^{-5/4}$.
So, since $a_n=o(n^{1/8})$, $(a_n+\|\pi\|_\infty)^2P(E_n^c)=o({n}^{-1})$.

Following result in Theorem \ref{th3} is provided by using Corollary \ref{coro2} and 
Corollary~\ref{coro3}.
\findemo

\section*{Acknoledgements}

I would like to thank F. Comte for her helpful suggestions throughout this work.

\bibliographystyle{elsart-harv}  
\bibliography{biblio}

\end{document}